# Ensemble Laplacian Biogeography-Based Sine Cosine Algorithm for Structural Engineering Design Optimization Problems

Vanita Garg[1], Kusum Deep[2], Khalid Abdulaziz Alnowibet[3], Ali Wagdy Mohamed[4,5], Mohammad Shokouhifar[6] and Frank Werner[7,*]

[1] School of Basic and Applied Sciences, Galgotias University, Greater Noida 201306, India
[2] Department of Mathematics, Indian Institute of Technology, Roorkee, Uttarakhand 247667, India
[3] Statistics and Operations Research Department, College of Science, King Saud University, Riyadh 11451, Saudi Arabia
[4] Operations Research Department, Faculty of Graduate Studies for Statistical Research, Cairo University, Giza 12613, Egypt
[5] Applied Science Research Center, Applied Science Private University, Amman 11931, Jordan
[6] Department of Electrical & Computer Engineering, Shahid Beheshti University, Tehran 1983969411, Iran
[7] Faculty of Mathematics, Otto-von-Guericke University, 39016 Magdeburg, Germany

**\* Correspondence:**
Frank Werner
Faculty of Mathematics, Otto-von-Guericke University, 39016 Magdeburg, Germany
Email: frank.werner@ovgu.de
Tel: +49 391 67 52025
Fax: +49 391 67 41171.

**Abstract:** In this paper, an ensemble metaheuristic algorithm (denoted as LX-BBSCA) is introduced. It combines the strengths of Laplacian Biogeography-Based Optimization (LX-BBO) and the Sine Cosine Algorithm (SCA) to address structural engineering design optimization problems. Our primary objective is to mitigate the risk of getting stuck in local minima and accelerate the algorithm's convergence rate. We evaluate the proposed LX-BBSCA algorithm on a set of 23 benchmark functions, including both unimodal and multimodal problems of varying complexity and dimensions. Additionally, we apply LX-BBSCA to tackle five real-world structural engineering design problems, comparing the results with those obtained using other metaheuristics in terms of objective function values and convergence behavior. To ensure the statistical validity of our findings, we employ rigorous tests such as the t-test and the Wilcoxon rank test. The experimental outcomes consistently demonstrate that the ensemble LX-BBSCA algorithm outperforms not only the basic versions of BBO, SCA, and LX-BBO but also other state-of-the-art metaheuristic algorithms.

# 1. Introduction

Optimization is considered a globally accepted tool for solving many real-world problems. Traditional optimization techniques have a limited scope in this regard as methods require continuity, differentiability knowledge of the objective function and constraints. Generally, real-world engineering problems are discontinuous and non-differentiable, and consequently, metaheuristic algorithms have been widely used in the applied sciences area because of their applicability and ease of use for a wide range of engineering problems [1-5].

Nature-inspired algorithms are optimization methods that are based on natural phenomena. Some popular metaheuristics in this class include Genetic Algorithm (GA), Ant Colony Optimization (ACO), Particle Swarm Optimization (PSO), and Grey Wolf Optimizer (GWO). These algorithms are rapidly evolving and growing. Based on the No-Free-Lunch theorem [6], this field has a lot of potential for researchers. It claims that no single algorithm is capable of solving all types of optimization problems, and each algorithm has its own advantages as well as some limitations.

Generally, exploration and exploitation are two important issues in stochastic algorithms [7]. Finding a new solution in the neighborhood of the existing solution is called exploration while making a high change in the current position of the existing solutions to find a better fitness value is termed as exploitation. Some of the metaheuristic algorithms focus on exploration while many of them are better at exploitation. A balanced optimizer is one that provides a better exploration and extensive exploitation. To achieve this goal, the hybridization of two or more metaheuristic algorithms has proved to be a successful idea in the literature.

The Sine Cosine Algorithm (SCA) [8] is a nature-inspired algorithm designed for solving complex optimization problems. It draws inspiration from the trigonometric functions of sine and cosine, which play a central role in the algorithm's behavior. The SCA employs a population of candidate solutions represented as angles and uses the sine and cosine functions to update their positions iteratively. This approach facilitates both exploration and exploitation of the solution space. SCA is known for its simplicity, ease of implementation, and efficiency in finding optimal or near-optimal solutions for a wide range of optimization problems, making it a valuable tool in optimization research and practical applications.

Biogeography-Based Optimization (BBO) [9] is a nature-inspired optimization algorithm rooted in the study of species distribution across geographic regions. It operates by modeling the migration and speciation of solutions within a solution space, akin to the movement of species between habitats. Each solution is associated with a Habitat Suitability Index (HSI), which is based on Suitability Index Variables (SIV) such as rainfall, vegetation, and climate. The BBO employs migration and dispersal operators to facilitate the transfer of high-quality solutions to less explored regions while introducing randomness for exploration. The BBO is adaptable to various problem domains, accommodating both continuous and discrete optimization tasks. BBO has found applications in diverse fields, including engineering, machine learning, and data mining, due to its ability to efficiently explore complex solution spaces, making it a valuable tool for real-world optimization challenges.

Laplacian BBO (LX-BBO) [10] is a variant of the BBO algorithm that incorporates Laplacian-based strategies to enhance its optimization capabilities. LX-BBO maintains the key principles of BBO, including the migration and speciation of solutions inspired by the species distribution in different habitats. However, LX-BBO introduces Laplacian operators to control the information flow among solutions in a population. These operators promote the sharing of useful information among the solutions while preventing overcrowding in the search space. By integrating Laplacian techniques, LX-BBO aims to improve the algorithm's convergence speed, exploration, and exploitation abilities, making it a more effective tool for solving complex optimization problems in various domains.

Despite the remarkable advancements in metaheuristic algorithms, it remains a fundamental truth that no single algorithm can universally excel at solving all types of optimization problems. Recognizing this limitation, the concept of hybridizing the SCA with the LX-BBO emerges as a promising approach for addressing a broad spectrum of real-world problems. SCA exhibits an intrinsic strength in exploration, while LX-BBO has demonstrated superior exploitation capabilities. By merging these two algorithms, we introduce a hybridized variant known as LX-BBSCA, seeking to harness the complementary attributes of both.

The structure of this paper is as follows: In Section 2, we give a comprehensive review of the current literature, exploring the landscape of existing metaheuristic algorithms. Section 3 offers an in-depth description of LX-BBO and SCA and introduces our innovative ensemble LX-BBSCA algorithm. Then in Section 4, we meticulously project the outcomes of our experimentation, encompassing an extensive analysis of the results, comprehensive statistical evaluations, insightful convergence analyses, and a thorough examination of the computational complexity. Section 5 takes a closer look at the application of the LX-BBSCA algorithm to tackle five structural engineering challenges, providing a detailed account of the solutions obtained. Finally, in Section 6, we draw our conclusions, summarizing the key findings of this study, and outline avenues for future research directions.

## 2. Literature review

The existing literature provides a comprehensive overview of the widespread applicability and enhancements of the SCA in various optimization domains. The SCA, introduced by Mirajalili in 2016 [8], has garnered attention for its simplicity and effectiveness in solving global optimization problems, as validated by numerous researchers. Notable applications of SCA include its utilization by Tawahid et al. [11] to address multi-objective optimization problems and its amalgamation with differential evolution (DE) by Bureerat and Pholdee [12] for structural damage detection in engineering. Furthermore, Reddy et al. [13] leveraged SCA to resolve profit-based unit commitment (PBUC) issues in competitive electricity markets, introducing a heuristic binary SCA method. Gupta and Deep [14] enhanced SCA by integrating PSO to improve its ability to locate global optima. Selim et al. [15] proposed a technique based on SCA and chaos map theory for allocating multiple distributed generators in distribution networks, demonstrating its practicality.

In hybridization efforts, SCA has been combined with various other algorithms to enhance its capabilities. Hota et al. [16] hybridized SCA with the Extreme Learning Machine (ELM) for exchange rate prediction, while Sindhu et al. [17] proposed a strategy to improve SCA's accuracy. Zhu et al. [18] introduced the orthogonal learning strategy (OLSCA) to enhance SCA, and Gupta and Deep [19] introduced the HSCA variant to mitigate local optima entrapment. In engineering design applications, El-Shorbagy et al. [20] harnessed the power of hybridization in SCA-SSGA, and Chen et al. [21] improved SCA's performance through the Nelder-Mead simplex concept and opposition-based learning. Singh et al. [22] proposed the GWO-SCA algorithm, combining GWO and SCA effectively to solve benchmark functions and engineering fields.

Additionally, several hybrid approaches have been explored, including the Sugeno Fuzzy Artificial Bee Colony algorithm [23], combined Whale Optimization Algorithm (WOA) and Simulated Annealing (WOA-SA) [24], and Heuristic-WOA-Variable Neighborhood Search (H-WOA-VNS) [25]. Fan et al. [26] introduced the SCA-FOA method, Mahdad et al. (2017) used SCA to enhance power system security, and Kumar et al. [27] utilized hybrid Cauchy and Gaussian mutations. Moreover, Beururat et al. [28] developed an adaptive SCA which is integrated with differential evaluation (ASCA-DE) to solve test problems for structural damage detection. Garg et al. [29] applied SCA to extract bioactive compounds, while Banerjee and Garg [30] proposed an algorithm incorporating five mutation operators into SCA for optimization tasks.

LX-BBO was extended by Garg and Deep [10], introducing constrained linear optimization capabilities, and Gupta and Deep [19] showcased LX-BBO in portfolio optimization. Furthermore, Garg et al. [31] demonstrated the effectiveness of incorporating LX-BBO into the Salp Swarm Algorithm and the Teaching Learning Algorithm. Gupta and Deep [19] combined SCA and the ABC algorithm, while Issa et al. [32] presented ASCA-PSO, a variant combining SCA with particle swarm optimization. Attia et al. [33] introduced Levy's flight into SCA to enhance its local search capabilities.

The extensive studies highlight the versatility and continuous improvements in the SCA and its various hybridizations, paving the way for enhanced optimization techniques across diverse domains. Das et al. [34] applied SCA to tackle fixed-head short-term hydrothermal scheduling problems, while Gupta and Deep [19] introduced the crossover approach (ISCA) to enhance SCA's exploration and exploitation abilities. The combination of SSGA's exploitation merits with SCA's exploration capabilities was proposed by Turgut et al. [35] for shell and tube evaporator design. Finally, Elaziz et al. [36] adapted SCA to address complex nonlinear optimization problems.

One of the notable trends is the fusion of SCA with other metaheuristics to leverage the strengths of multiple algorithms, resulting in enhanced exploration and exploitation capabilities. It indicates a pragmatic approach to solving complex problems by capitalizing on the complementary attributes of different algorithms. Furthermore, the incorporation of various mutation strategies, including Cauchy, Gaussian, and Levy mutations, has been demonstrated to augment SCA's ability to navigate intricate solution spaces efficiently. These mutation operators introduce diversity and adaptability into the optimization process, potentially mitigating convergence to local optima.

The literature underscores the importance of rigorous testing and validation through benchmark functions and practical engineering applications. By subjecting SCA and its variants to diverse testing scenarios, researchers have provided empirical evidence of their effectiveness in solving complex problems. In conclusion, the diverse body of literature demonstrates the ever-evolving landscape of SCA and its hybridizations, reinforcing its status as a robust and versatile optimization tool. The reviewed literature highlights the potential directions of SCA such as exploring hybridization with emerging metaheuristic algorithms and conducting more extensive real-world applications to further validate its practical utility.

## 3. Proposed LX-BBSCA algorithm

*3.1. Laplacian biogeography-based optimization (LX-BBO)*

BBO [9] is based on the island biogeography theory, which is concerned with the migration, speciation, and extinction of species within a habitat. It only considers the better living conditions for measuring the HSI. Personal interests, individual choices, or societal norms, are not taken into consideration. A habitat or living place is considered a solution. The solution's fitness value is the number of species in the habitat. More species in a habitat is a notion of a good solution as this habitat has better living conditions and thus attracts species to move into the place. This phenomenon is termed a migration operator, as the better habitat is preferred over a worse solution, i.e., a habitat that has a smaller number of species count. The drastic changes in weather (e.g., tsunami, earthquake, or volcano eruption) is termed mutation operator in BBO. The details of these operators can be given in the following.

3.1.1. Migration operator

Based on emigration and immigration rates, the solutions are being replaced with a better solution. High emigration rates indicate a better solution and a high immigration rate gives a worse solution. This operator is an important part of the exploitation of the search space. A balanced and robust algorithm has a strong exploitation so as to find solutions in the neighborhood of the existing solution. In classical BBO, the migration operator is lagging behind to fully exploit the search space. It excels at exploitation but falls short of exploration. The BBO has undergone numerous modifications and hybridizations since its inception in order to improve its performance. There have been numerous efforts to raise the standard of BBO's exploitation. In an effort to better leverage the search space, Garg and Deep [37] collected various types of upgraded BBO utilizing different mutation and migration operators. A real-coded genetic algorithm's LX-crossover by Garg and Deep [10] produced a remarkable LX-BBO performance.

3.1.2. Laplacian migration operator

The restriction of the migration operator in the original BBO has been solved by the addition of the LX crossover from the real coded GA. The LX-BBO results have properly asserted that the program performed superbly. The LX-BBO migration operator selects two habitats ($x_1$ and $x_2$) using emigration and immigration rates. Then, these habitats will be combined via Laplace crossover to produce new habitats. The two new habitats, denoted as $y_1$ and $y_2$, are produced as follows:

$$y_1^i = x_1^i + \beta(x_1^i - x_2^i) \tag{1}$$

$$y_2^i = x_2^i + \beta(x_1^i - x_2^i) \tag{2}$$

where $\beta$ is a random number that follows the Laplace distribution and is generated as follows:

$$\beta = \begin{cases} a - b\log(u) & if\ u \leq 1/2 \\ a + b\log(u) & if\ u > 1/2 \end{cases} \qquad (3)$$

where $u$ is a uniform random number within [0 1], $b > 0$ is called the scale parameter, and $a \in R$ is a location parameter. The habitats $y_1$ and $y_2$ produce a new habitat $z$ for the next generations. It is done as follows:

$$z = \gamma y_1^i + (1-\gamma)y_2^i \qquad (4)$$

$$\gamma = \gamma_{min} + (\gamma_{max} - \gamma_{min})^k \qquad (5)$$

where $\gamma_{min}$ and $\gamma_{max}$ are in the range [0, 1] and define the minimum and maximum bounds of $\beta$, $i$ is the counter of the number of generations, and $k$ is a user-specific parameter. After each generation, an Elitism mechanism is applied to select two habitats. It preserves the two best environments to replace the two worst solutions. The pseudo-code of LX-BBO is shown in Algorithm 1.

### 3.1.3. Mutation operator

Mutation in classical BBO is the tool for exploring new solutions. The newly generated solution by the mutation operator is based on the mutation rate and these rates are determined by a differential equation. A rigorous experimental study on the effects of different mutations on the performance of LX-BBO can be found in [38]. However, the results are almost similar for all mutation operators, and no outstanding breakthrough was established. This proves that the mutation operator needs to be efficient for the LX-BBO algorithm.

---

Algorithm 1. LX-BBO algorithm.

**Begin**
Generation of a random population H={$H_1, H_2, H_3, \ldots, H_n$}
Determining the corresponding HSI values
Sort the population by the HSI values from best to worst
while ($t \leq Max\_iteration$)
    for each habitat
        for each SIV
            Calculation of the emigration and immigration rates for each habitat
            Choosing $H_i$ based on the immigration rate
            if $H_i$ is chosen $(i = 1,2,3,\ldots,n)$
                $H_i^1(SIV) = H_i(SIV) + \beta(H_i(SIV) - H_j(SIV))$
                $H_i^2(SIV) = H_j(SIV) + \beta(H_i(SIV) - H_j(SIV))$
                Updating $H_i(SIV) \leftarrow \gamma H_i^1(SIV) + (1-\gamma)H_i^2(SIV)$
                $\gamma = \gamma_{min} + (\gamma_{max} - \gamma_{min})^k$
                if $H_j$ is chosen $(j = 1,2,3,\ldots,n)$
                      Randomly selection of an SIV from $H_i$
                      Replacing a randomly SIV in $H_i$ with one from $H_j$
                end if

---

```
            end if
            Choosing a SIV in $H_i$ based on the mutation rate $\mu$
            if $H_j(SIV)$ is chosen
                Replacing a randomly selected SIV
            end if
        end for
    end for
    Calculation of the modified habitats HSI value again
    Sort the population by the HSI values from best to worst
    Replacing the worst solution using Elitism strategy
    Sort the population by the HSI values from best to worst
    t=t+1;
end while
End
```

## 3.2. Sine cosine algorithm (SCA)

The SCA was inspired by the sine and cosine functions from mathematics and is designed to solve optimization problems. It falls under the category of metaheuristic algorithms, which are used to find approximate solutions to optimization and search problems. This algorithm simulates the movement of sine and cosine waves to guide the search for optimal solutions in a given problem space. It is a population-based algorithm, where each potential solution is represented as an individual in the population. The algorithm iteratively updates the positions of these individuals in search of better solutions over time. Next, we present the basic outline of how the SCA works:

1. *Initialization*: Initialize a random population of candidate solutions within the search space.
2. *Evaluation*: Evaluate the fitness of the candidate solutions based on the objective function. The fitness value indicates how well each solution performs in terms of the objective.
3. *Update positions*: Update the position of each candidate solution using the following equations, which are inspired by the sine and cosine functions:

$$X_i^{t+1} = \begin{cases} X_i^t + r_1 sin(r_2)|r_3 P_i^t - X_i^t| & r_4 \leq 0.5 \\ X_i^t + r_1 cos(r_2)|r_3 P_i^t - X_i^t| & r_4 > 0.5 \end{cases} \quad (6)$$

where $X_i^t$ is the position of solution *i* in the current iteration, $X_i^{t+1}$ is the position of the same solution in the next iteration, $P_i^t$ is the destination position of the global optimal particle, $r_1$ linearly decreases from 2 to 0, and $r_2$, $r_3$ and $r_4$ are uniformly distributed random numbers: $r_2$ and $r_3 \in [0,2]$; $r_4 \in [0,1]$.

4. *Boundary handling*: It ensures that the updated positions remain within the feasible region of the problem by applying suitable boundary handling mechanisms.
5. *Update the best solution*: It updates the best solution found so far based on the fitness values of the new candidate solutions.
6. *Termination*: The algorithm terminates when a certain stop criterion is met, e.g., reaching a pre-defined number of iterations or achieving a satisfactory solution.

SCA is relatively simple to implement and has shown promising performance on different problems. However, its performance can vary depending on the problem's characteristics and parameter settings. As with any optimization algorithm, proper parameter tuning and problem-specific adaptations are crucial for achieving good results. This method considerably assists the search process in avoiding local optima, resulting in faster search times. Its exploration abilities are superior to its exploitation abilities. The pseudo-code of SCA is shown in Algorithm 2. In this pseudocode, $H_i^t$ denotes habitant $i$ in iteration $t$.

---

Algorithm 2. Classical CSA algorithm.

**Begin**
Initialization of the solutions $H_i (i = 1,2,3, \dots, n)$
Assessment of the fitness value of each solution
Selecting the best solution as $Gbest$
while ($t \leq Max\_iteration$)
    for $i = 1:n$
        Updating $r_1$ to $r_4$: $r_1$: linearly decreasing from 2 to 0, $r_2$ and $r_3 \in [0,2]$, $r_4 \in [0,1]$
        Updating the position of solution $i$ using Eq. (6)
        Assessment of the fitness value for solution $i$
        Checking the solution whether it goes beyond the search space; if yes, bring it back
        if $f(H_i^{t+1}) < f(H_i^t)$
            Updating the global best solution $Gbest$
        end if
    end for
    *t=t*+1;
end while
return $Gbest$
**End**

---

*3.3. Proposed hybrid version of LX-BBO and SCA (LX-BBSCA)*

Algorithm LX-BBO is good enough at altering the current population given to its migration operator, but it takes a long time to explore the entire search space. Therefore, a combined algorithm based on LX-BBO and SCA is designed to choose the most relevant feature set from the original feature. These algorithms are combined into a hybrid version, denoted as LX-BBSCA, where the mutation operator of LX-BBO is replaced with that of SCA. The pseudo-code of the LX-BBSCA is given in Algorithm 3.

Algorithm 3. Proposed LX-BBSCA algorithm.

**Begin**
Initialization of the solutions $H_i (i = 1,2,3, \ldots, n)$
Assessment of the fitness value of each solution
Selecting the best solution as $Gbest$
while ($t \leq Max\_iteration$)
    SCA:
    for $i = 1:n$
        Updating $r_1$ to $r_4$: $r_1$: linearly decreasing from 2 to 0, $r_2$ and $r_3 \in [0,2]$, $r_4 \in [0,1]$
        Updating the position of solution $i$ using Eq. (6)
        Assessment of the fitness value for solution $i$
        Checking the solution whether it goes beyond the search space; if yes, bring it back
        if $f(H_i^{t+1}) < f(H_i^t)$
            Updating the global best solution $Gbest$
        end if
    end for
    LX-BBSCA:
    Calculate the immigration rate
    Select habitant $H_i$
    for $i = 1:n$
        for $j = 1:D$
            if rand<$\lambda_k$
                $H_i^1(SIV) = H_i(SIV) + \beta(H_i(SIV) - H_j(SIV))$
                $H_i^2(SIV) = H_j(SIV) + \beta(H_i(SIV) - H_j(SIV))$
                Updating $H_i(SIV) \leftarrow \gamma H_i^1(SIV) + (1 - \gamma) H_i^2(SIV)$
                $\gamma = \gamma_{min} + (\gamma_{max} - \gamma_{min})^k$
            end if
            if rand<0.5 then
                Replacing the SIV of $H_i$ with a randomly selected SIV
            end if
        end for
    end for
    for each habitat
        Recompute the HSI values
        Checking the solution whether it goes beyond the search space; if yes, bring it back
        if $f(H_i^{t+1}) < f(H_i^t)$
            Updating the global best solution $Gbest$
        end if
    end for
    *t=t+*1;
end while
Return $Gbest$
**End**

## 4. Validation of LX-BBSCA for benchmark functions

To validate the performance of the LX-BBSCA algorithm, it is used to solve 23 popular benchmark functions, as listed in Table 1. The problem set consists of 7 unimodal test functions (F1–F7) and 16 multimodal test functions (F8–F23).

**Table 1.** Benchmark functions.

| Function | Function name | Dim | Range | $f_{min}$ |
|---|---|---|---|---|
| $f_1(x) = \sum_{i=1}^{n} x_i^2$ | sphere | 30 | [-100,100] | 0 |
| $f_2(x) = \sum_{i=1}^{n} |x_i| + \prod_{i=1}^{n} |x_i|$ | Schwefel's 2.22 | 30 | [-10,10] | 0 |
| $f_3(x) = \sum_{i=1}^{n} x_i (\sum_{i=1}^{n} x_i)^2$ | Schwefel's 1.20 | 30 | [-100,100] | 0 |
| $f_4(x) = \max\{|x_i|, 1 \le i \ge n$ | Schwefel's 2.21 | 30 | [-100,100] | 0 |
| $f_5(x) = \sum_{i=1}^{n-1}[100(x_{i+1} - x_i^2)^2 + (x_i - 1)^2]$ | Rosenbrock | 30 | [-30,30] | 0 |
| $f_6(x) = \sum_{i=1}^{n}(|x_i + 0.5|)^2$ | step | 30 | [-100,100] | 0 |
| $f_7(x) = \sum_{i=1}^{n} ix_i^4 + random[0,1)$ | Quartic Noise | 30 | [-1.28,1.28] | 0 |
| $f_8(x) = \sum_{i=1}^{n} -x_i \sin(\sqrt{|x_i|})$ | Schwefel | 10 | [-500,500] | -418*5 |
| $f_9(x) = \sum_{i=1}^{n}[x_i^2 - 10\cos(2\pi x_i) + 10]$ | Rastrigin | 30 | [-5.12,5.12] | 0 |
| $f_{10}(x) = -20exp\left(-0.2\sqrt{\frac{1}{n}\sum_{i=1}^{n} x_i^2}\right) -$ $exp(\frac{1}{n}\sum_{i=1}^{n} \cos(2\pi x_i)) + 2 + e$ | Ackley | 30 | [-32,32] | 0 |
| $f_{11}(x) = \frac{1}{4000}\sum_{i=1}^{n} x_i^2 - \prod_{i=1}^{n} \cos\left(\frac{x_i}{\sqrt{i}}\right) + 1$ | Griewank | 30 | [-600,600] | 0 |
| $f_{12}(x) = \frac{\pi}{n}\{10\sin(\pi y_i) + \sum_{i=1}^{n-1}(y_i - 1)^2[1 +$ $10sin^2(\pi y_{i+1})] + (y_n - 1)^2\} + \sum_{i=1}^{n} u(x_i, 10, 100, 4)$ | Penalized 1 | 30 | [-50,50] | 0 |
| $f_{13}(x) = 0.1\{sin^2(3\pi x_1) + \sum_{i=1}^{n}(x_i - 1)^2[1 + sin^2(3\pi x_i + 1)] + (x_n - 1)^2[1 + sin^2(2\pi x_n)]\} + \sum_{i=1}^{n} u(x_i, 5, 100, 4)$ | Michalewiez | 30 | [-50,50] | 0 |
| $f_{14}(x) = \left(\frac{1}{500} + \sum_{j=1}^{25}\frac{1}{j+\sum_{i=1}^{2}(x_i-a_{ij})^6}\right)^{-1}$ | Shekel's | 2 | [-65,65] | 1 |
| $f_{15}(x) = \sum_{i=1}^{11}[a_i - \frac{x_1(a_i^2+b_i x_2)}{b_i^2+b_i x_3+x_4}]^2$ | Foxholes | 2 | [-5,5] | 0.00030 |
| $f_{16}(x) = -4x_1^2 - 2.1x_1^4 + \frac{1}{3}x_1^6 + x_1 x_2 - 4x_2^4 + 4x_2^4$ | Six-hump camel | 3 | [-5,5] | 0.398 |
| $f_{17}(x) = (x_2 - \frac{5.1}{4\pi^2}x_1^2 + \frac{5}{\pi}x_1 - 6)^2 + 10(1 - \frac{1}{8\pi})cosx_1 + 10$ | Goldstein price | 6 | [-5,5] | 3 |
| $f_{18}(x) = [1 + (x_1 + x_2 + 1)^2(19 - 14x_1 + 3x_1^2 - 14x_2 + 6x_1 x_2 + 3x_2^2)] * [30 + (2x_1 - 3x_2)^2]$ | Hartman's family-2 | 4 | [0,10] | -3.86 |
| $f_{19}(x) = -c_i exp(-\sum_{i=1}^{6} a_{ij}(x_j - P_{IJ})^2)$ | sphere | 3 | [1,3] | -3.32 |
| $f_{20}(x) = -\sum_{i=1}^{4} c_i \exp(-\sum_{i=1}^{3} a_{ij}(x_j - P_{IJ})^2)$ | Hartmann 6-D | 6 | [0,1] | -10.15 |
| $f_{21}(x) = -\sum_{i=1}^{5}[(X - a_i)(X - a_i)^T + C^i]^{-1}$ | Shekel's | 4 | [0,10] | -10.15 |
| $f_{22}(x) = -\sum_{i=1}^{7}[(X - a_i)(X - a_i)^T + C^i]^{-1}$ | -- | 4 | [0,10] | -10.40 |
| $f_{23}(x) = -\sum_{i=1}^{10}[(X - a_i)(X - a_i)^T + C^i]^{-1}$ | -- | 4 | [0,10] | -10.53 |

In the following, the results of the proposed LX-BBSCA algorithm are compared to BBO and LX-BBO. The statistical parameters include the values for the objective function's average, the standard deviation (Std), the median, the minimum and maximum values, as well as the minimum and maximum for each run. The validity of the algorithms was evaluated using the Std and average. Furthermore, the minimum and maximum values of the objective function represent the best possible solutions to the given problem in terms of iterations. In Table 2, the calculations are computed for 30 independent trials and 1000 function evaluations.

**Table 2.** Objective function values obtained by LX-BBO, BBO, and LX-BBSCA.

| Function | Algorithm | Min | Max | Std | Average | Median |
|---|---|---|---|---|---|---|
| 1 | LXBBSCA | **1.47E-08** | **0.003555** | **0.001019** | **0.0005** | **3.15E-05** |
|   | BBO | 1777.031 | 5077.365 | 716.4427 | 3868.813 | 3941.62 |
|   | LX-BBO | 1494.216 | 5044.481 | 865.3063 | 3672.611 | 3714.437 |
| 2 | LXBBSCA | **2.03E-06** | **0.004806** | **0.000867** | **0.000304** | **7.91E-05** |
|   | BBO | 10.79175 | 21.72683 | 2.374523 | 17.00374 | 17.2588 |
|   | LX-BBO | 11.32607 | 19.53608 | 2.057063 | 16.08974 | 16.2298 |
| 3 | LXBBSCA | **0.006513** | **13.14504** | **3.474535** | **1.89781** | **0.278841** |
|   | BBO | 1934.879 | 5951.809 | 927.4112 | 3758.93 | 3920.839 |
|   | LX-BBO | 2049.36 | 5171.657 | 816.2047 | 3564.492 | 3472.326 |
| 4 | LXBBSCA | **0.006297** | **0.446159** | **0.12764** | **0.113329** | **0.043289** |
|   | BBO | 22.23739 | 42.85724 | 4.189265 | 34.69546 | 34.97368 |
|   | LX-BBO | 19.27473 | 39.92089 | 4.00352 | 33.3647 | 33.24844 |
| 5 | LXBBSCA | **7.505546** | **706.7561** | **134.2925** | **40.99334** | **8.27666** |
|   | BBO | 437754 | 3670960 | 763242.4 | 1595315 | 1487924 |
|   | LX-BBO | 672237.7 | 2843717 | 561336.2 | 1672675 | 1735165 |
| 6 | LXBBSCA | **0.288575** | **0.979075** | **0.153504** | **0.703137** | **0.731203** |
|   | BBO | 1932.084 | 6250.072 | 845.4571 | 3853.714 | 3971.152 |
|   | LX-BBO | 987.377 | 5581.525 | 888.0294 | 3859.914 | 3877.024 |
| 7 | LXBBSCA | **0.001156** | **0.027382** | **0.006862** | **0.007034** | **0.003887** |
|   | BBO | 0.26968 | 1.230022 | 0.230808 | 0.758662 | 0.709212 |
|   | LX-BBO | 0.226661 | 1.362734 | 0.30781 | 0.809159 | 0.840667 |
| 8 | LXBBSCA | **-2827.44** | **-2036.08** | **166.5375** | **-2244.35** | **-2216.65** |
|   | BBO | -2829.83 | -2121.53 | 181.1336 | -2456.98 | -2445.76 |
|   | LX-BBO | -2875.34 | -2225.78 | 131.6692 | -2432.53 | -2413.6 |
| 9 | LXBBSCA | **2.33E-06** | **35.86702** | **10.85302** | **5.389507** | **0.14953** |
|   | BBO | 46.81987 | 78.79715 | 6.388578 | 65.23132 | 65.17603 |
|   | LX-BBO | 40.7377 | 74.05955 | 6.887215 | 62.28468 | 63.68432 |
| 10 | LXBBSCA | **6.82E-05** | **0.056023** | **0.010371** | **0.004476** | **0.001321** |
|   | BBO | 13.59016 | 17.02045 | 0.828409 | 15.85072 | 16.04995 |
|   | LX-BBO | 13.5316 | 17.03005 | 0.923063 | 15.67928 | 15.87936 |
| 11 | LXBBSCA | **0.00015** | **0.863571** | **0.231311** | **0.330496** | **0.313525** |
|   | BBO | 22.21162 | 48.54429 | 7.655131 | 35.48117 | 35.48288 |
|   | LX-BBO | 11.99046 | 48.13176 | 7.368682 | 36.91168 | 37.60413 |
| 12 | LXBBSCA | **0.066566** | **0.350893** | **0.070983** | **0.182471** | **0.179451** |
|   | BBO | 42333.56 | 2511312 | 685507 | 686604.5 | 373467.9 |
|  | LX-BBO | 266.1032 | 2036870 | 468695.3 | 488556.4 | 353698.8 |
| 13 | LXBBSCA | **0.199118** | **0.547923** | **0.085934** | **0.408174** | **0.401248** |
|  | BBO | 697435.3 | 12295841 | 2913431 | 4140497 | 3615374 |
|  | LX-BBO | 475863.4 | 9962896 | 2547311 | 4441777 | 4136487 |
| 14 | LXBBSCA | **0.998004** | **1.009151** | **0.002023** | **0.998519** | **0.998021** |
|  | BBO | 0.998006 | 1.027821 | 0.006486 | 1.000994 | 0.998702 |
|  | LX-BBO | 0.998004 | 1.176247 | 0.041217 | 1.019683 | 1.001262 |
| 15 | LXBBSCA | **0.000522** | **0.001586** | **0.000381** | **0.001096** | **0.001103** |
|  | BBO | 0.000736 | 0.004441 | 0.000861 | 0.002095 | 0.00193 |
|  | LX-BBO | 0.000443 | 0.005199 | 0.001079 | 0.002494 | 0.002101 |
| 16 | LXBBSCA | **-1.03162** | **-1.03129** | **8.75E-05** | **-1.03155** | **-1.03158** |
|  | BBO | -1.03152 | -1.01417 | 0.004388 | -1.02788 | -1.02916 |
|  | LX-BBO | -1.03136 | -1.02386 | 0.002109 | -1.02922 | -1.02999 |
| 17 | LXBBSCA | **0.397887** | **0.407469** | **0.002234** | **0.399293** | **0.39849** |
|  | BBO | 0.397969 | 0.402253 | 0.001318 | 0.399398 | 0.398989 |
|  | LX-BBO | 0.397943 | 0.405405 | 0.002047 | 0.399853 | 0.39913 |
| 18 | LXBBSCA | **3.000001** | **3.001711** | **0.000314** | **3.000123** | **3.000031** |
|  | BBO | 3.000651 | 3.301352 | 0.070737 | 3.062178 | 3.039374 |
|  | LX-BBO | 3.001377 | 3.208721 | 0.060075 | 3.055511 | 3.026186 |
| 19 | LXBBSCA | **-3.86153** | **-3.85282** | **0.001996** | **-3.85482** | **-3.85442** |
|  | BBO | -3.8626 | -3.85139 | 0.002833 | -3.8588 | -3.85956 |
|  | LX-BBO | -3.86254 | -3.84367 | 0.004217 | -3.85762 | -3.8583 |
| 20 | LXBBSCA | **-3.12678** | **-1.68596** | **0.347211** | **-2.88178** | **-3.00666** |
|  | BBO | -3.22558 | -2.95942 | 0.061642 | -3.0718 | -3.05945 |
|  | LX-BBO | -3.19627 | -2.91769 | 0.056134 | -3.08601 | -3.08811 |
| 21 | LXBBSCA | **-5.85902** | **-0.87934** | **1.738807** | **-3.60057** | **-4.69579** |
|  | BBO | -8.40559 | -2.53527 | 1.188802 | -4.05523 | -3.75194 |
|  | LX-BBO | -7.07924 | -2.17446 | 1.127562 | -3.73127 | -3.67004 |
| 22 | LXBBSCA | **-8.8799** | **-0.90652** | **1.968125** | **-4.25467** | **-4.60085** |
|  | BBO | -5.87171 | -2.48824 | 0.819649 | -3.90592 | -3.90487 |
|  | LX-BBO | -8.17891 | -2.59073 | 1.288829 | -4.16432 | -3.82846 |
| 23 | LXBBSCA | **-6.73001** | **-0.94429** | **1.530706** | **-4.24441** | **-4.38713** |
|  | BBO | -7.76979 | -2.86094 | 1.138833 | -4.29905 | -4.01435 |
|  | LX-BBO | -6.95564 | -2.57408 | 1.23627 | -4.39294 | -4.20171 |

*4.1. Unimodal test functions*

If an unimodal test function is monotonically increasing and monotonically decreasing, respectively, then the function has one global minimum or maximum value. The unimodal functions (F1–F7) have one global solution. So, they are taken to investigate the exploitation capacity. As a result, they are used to assess the convergence rate of the heuristic optimization algorithms. In terms of the mean, the Std median, the minimum and maximum values of the result.

To solve an unimodal problem, the hybrid method LX-BBSCA can find the most feasible and comprehensive solutions. In all problems, Algorithm LX-BBSCA outperforms the other algorithms. As a result, it has a strong exploitation ability and great performance when it locates the overall solution of the unimodal benchmark function. In Table 2, the objective function values are analyzed over 30 different runs, and their mean, mode, Std, best values, and worst values are reported.

The results in the first seven rows give the results for unimodal functions. Clearly, LXBBSCA has obtained better results than its counterparts. The results indicate that LXBBSCA has better exploitation capability than LXBBO and BBO.

*4.2. Multimodal test functions*

The multimodal functions (F8–F13) have a lot of local optimal points. These functions are typically tested to assess the exploration ability of the different methods, as a lot of local optima values raises the possibility of stagnation of an algorithm. Multimodal functions with numerous local minima are notoriously difficult to optimize. The proposed variant's superiority and ability on multiple dimensions have been demonstrated using the Std, the average, and the minimum and maximum objective function values. The final results are significantly more relevant for these functions since they represent the algorithm's ability to escape undesirable local optima and obtain a good solution near the global optimum.

All LX-BBSCA techniques outperformed BBO, LX-BBO, and SCA in the final finding for all six functions (F8–F13). Furthermore, the test results reveal that the novel hybrid method's strong exploration capability allows us to extensively search the space and find viable search areas. In Table 2, the results are demonstrated for multimodal functions in rows 8-13. The analysis of these problems clearly indicates the better performance of LXBBCA as the objective function values of LX-BBSCA are smaller than those obtained by LX-BBO and BBO. Since the problems are minimization ones, the results are encouraging in order to claim that LX BBCA has a better exploration strategy than its counterparts. Multimodal problems need a better extensive search in the search space.

The test problems (F14–F23) are multimodal functions that are typically used to evaluate the search algorithm's ability to balance exploration and exploitation. The main difference between functions F8–F13 and functions with a few local minima, i.e., F14–F23, is that functions F14–F23 are simpler than functions F8–F13 due to their low dimensionality and fewer local minima. Algorithms BBO and LX-BBO were significantly outperformed by LX-BBSCA.

These findings illustrate that the proposed LX-BBSCA algorithm outperforms the traditional method in terms of solution quality. The new hybrid approach LX-BSCA prevents the search process from converging prematurely to a local optimal point and allows exploration of the search path. In Table 2, rows 14-23, the results are given for multimodal functions with varying dimensions. Testing the algorithm on varying dimensionality problems is the key to finding the usability and applicability of the algorithm to complex real-world problems. The results in Table 2 clearly show the better performance of the proposed method over its basic versions.

*4.3. Statistical analysis*

To check the performance of the algorithms, the t-test and the Wilcoxon rank test are performed, and compared the significance of the algorithms LX-BBSCA and LX-BBO.

4.3.1.  t-test analysis

A paired t-test, also referred to as a paired-samples t-test or a dependent-samples t-test, is a statistical test employed to assess whether there exists a significant disparity between the means of two interconnected groups. This test is particularly useful when one has two sets of observations that are paired in some way, such as before-and-after measurements on the same subjects or matched data points. It is appropriate when the data are assumed to be approximately normally distributed and the differences between the paired observations are approximately normally distributed as well.

Table 3 shows the results obtained after using the t-test for the pairs of Algorithms LX-BBSCA and LX-BBO. The p-values are significant because they provide a common language for the test results. A p-value is a number ranging from 0 to 1 which can be used to determine the statistical significance of a set of data.

1. If the p-value is larger than 0.05, the difference described is not significant.
2. If the p-value is less than 0.05 but greater than 0.001, the difference is significant.
3. If the p-value is less than 0.001, the difference is very significant.

The result analysis of the hybrid algorithm LX-BBSCA and LX-BBO are shown in Table 3, in which, "a+", "a", and "b" mean the difference is "very significant", "significant", and "no significant", respectively.

**Table 3.** t-test results for LX-BBSCA and LX-BBO.

| Function | Mean | Std Deviation | Std Error Mean | 95% Confidence Interval | | Significance | |
|---|---|---|---|---|---|---|---|
| | | | | Lower | Upper | p | Conclusion |
| F1 | 3672.61 | 865.306 | 157.56 | 3995.72 | 3349.5 | 0 | a+ |
| F2 | 16.0892 | 2.057 | 0.37558 | 16.8574 | 15.3211 | 0 | a+ |
| F3 | 3562.59 | 816.499 | 149.071 | 3867.48 | 3257.71 | 0 | a+ |
| F4 | 33.2513 | 4.02652 | 0.73514 | 34.7549 | 31.7478 | 0 | a+ |
| F5 | 1672634 | 561349 | 102488 | 1882245 | 1463023 | 0 | a+ |
| F6 | 3859.21 | 888.051 | 162.135 | 4190.81 | 3527.61 | 0 | a+ |
| F7 | 8E+07 | 0.30838 | 0.0563 | 9172777 | 6.9E+07 | 0 | a+ |
| F8 | 188.181 | 206.407 | 37.6847 | 11.1075 | 265.255 | 0 | a+ |
| F9 | 56.895 | 11.0961 | 2.02585 | 61.03 | 52.7518 | 0 | a+ |
| F10 | 15.6748 | 0.92076 | 0.16811 | 16.0186 | 15.3309 | 0 | a+ |
| F11 | -36.581 | 7.38748 | 1.34876 | 39.3397 | 33.8227 | 0 | a+ |
| F12 | -488556 | 468695 | 85571.7 | 663570 | 313542 | 0 | a+ |
| F13 | 4441776 | 2547311 | 465073 | 5392958 | 3490595 | 0 | a+ |
| F14 | 0.0211647 | 0.04141 | 0.00756 | -4E+06 | 0.005701113 | 0.009 | a |
| F15 | 13988 | 0.00108 | 0.0002 | -0.0018 | -0.001 | 0 | a+ |
| F16 | 0.00233 | 0.00213 | 0.00039 | 0.00313 | 0.00154 | 0 | a+ |
| F17 | 0.00055987 | 0.00317 | 0.00058 | 17442 | 0.00062 | 0.034 | a |
| F18 | 0.05539 | 0.05998 | 0.01095 | 0.07779 | 0.03299001 | 0 | a+ |
| F19 | 0.00279 | 0.00447 | 0.00082 | 0.00112 | 0.00446 | 0.002 | a |
| F20 | 200620 | 0.35852 | 0.06546 | 0.06675 | 0.3345 | 0.005 | a |
| F21 | 0.15514 | 2.13066 | 0.38901 | 6.4E+09 | 0.95075 | 0.693 | b |
| F22 | 0.035 | 2.2221 | 0.40571 | -0.8648 | 0.79477 | 0.932 | b |
| F23 | 0.14853 | 2.1494 | 0.39243 | 0.65407 | 0.95114 | 0.708 | b |

### 4.3.2. Wilcoxon rank sum test

Wilcoxon's rank-sum test (WRST) is utilized to provide a more statistical evaluation of the suggested approach. This is a non-parametric test that determines whether LX-BBSCA and LX-BBO have a significant difference. For a statistically significant correlation, the Wilcoxon test was used since the median value is more essential in describing the algorithm's efficiency, and this test may be applied without knowing the distribution of the data set. The statistical comparison between LX-BBO and the proposed LX-BBSCA algorithm is done in Table 4 at a significance level of 5%. The acquired p-values corresponding to all test problems are also presented in this table. The "+" sign is used to show that the LX-BBSCA is statistically significantly better than LX-BBO, and the "-" sign is considered to mention that LX-BBO is statistically significantly better than LX-BBSCA. Table 4 shows that out of 23 benchmark test problems, the proposed LX-BBSCA algorithm significantly outperforms LX-BBO in 17 test functions, while LX-BBO significantly outperforms LX-BBSCA in the remaining seven problems.

**Table 4.** Wilcoxon test for LX-BBSCA and LX-BBO with 5% of significance.

| Function | Wilcoxon's Z value | p | Conclusion |
| --- | --- | --- | --- |
| F1 | -4.782 | 0.000 | + |
| F2 | -4.782 | 0.000 | + |
| F3 | -4.782 | 0.000 | + |
| F4 | -4.782 | 0.000 | + |
| F5 | -4.782 | 0.000 | + |
| F6 | -4.782 | 0.000 | + |
| F7 | -3.815 | 0.000 | + |
| F8 | -4.78 | 0.000 | + |
| F8 | -4.782 | 0.000 | + |
| F9 | -4.782 | 0.000 | + |
| F10 | -4.782 | 0.000 | + |
| F11 | -4.782 | 0.000 | + |
| F12 | -4.782 | 0.000 | + |
| F13 | -3.815 | 0.000 | + |
| F14 | -3.78 | 0.000 | + |
| F15 | -3.887 | 0.000 | + |
| F16 | -0.234 | 0.000 | + |
| F17 | 0.256 | 0.245 | - |
| F18 | 0.345 | 0.000 | + |
| F19 | 0.332 | 0.003 | - |
| F20 | 0.3075 | 0.002 | - |
| F21 | -0.237 | 0.813 | - |
| F22 | 0.319 | 0.943 | - |
| F23 | 0.319 | 0.750 | - |

*4.4. Computational complexity*

The computational complexity of an algorithm is a criterion for determining its robustness with the least amount of effort. The computational complexity of the LX-BBO and LX-BBSCA algorithms can be described as follows. Initially, both algorithms consider the same population size $N$. Both algorithms have the same complexity when it comes to initializing input parameters such as population size, number of iterations, problem size, and so on. After a while loop has been completed, the main mathematical measures for both algorithms begin. Because each solution is updated in Algorithms LX-BBO and LX-BBSCA, this computational complexity is $O(N)$. If the algorithm fails to achieve the desired accuracy, it will be continued for $T$ iterations. Moreover, an objective function is used to determine each solution. As a result, the computational complexity in this step is $O(N)$. In the worst-case scenario, the overall complexity can be defined as $O(T*N)$.

## 5. Applying LX-BBSCA to structural engineering design problems

This section evaluates Algorithms LX-BBO and BBO with the proposed Algorithm LX-BBSCA for three structural engineering problems and two reliability problems.

*5.1. Gear train problem*

Analyzing the connections between several gears in a mechanical system is the gear train problem. Between spinning shafts, motion, and power are transferred using gears. They come in a variety of shapes and sizes, and depending on how they are arranged in a gear train, they can offer varying mechanical benefits, speed ratios, and torque variations. Unconstrained problems include the Gear Train Problem from 1990. A gear train is defined as a set of gears that consists of two pairs of gears and has an overall reduction ratio that as closely as possible approaches a certain value while keeping in mind that the total number of teeth in all four gears must not exceed a predetermined limit. The problem can be mathematically stated as follows:

$$\text{minimize } f_1(x) = \left(\frac{1}{6.931} - \frac{x_1 x_3}{x_2 x_4}\right)^2, \quad x = (x_1, x_2, x_3, x_4) \quad (7)$$

$$s.t. \quad 12 \leq (x_1, x_2, x_3, x_4) \leq 60 \quad (8)$$

where the variables $x_1$, $x_2$, $x_3$, and $x_4$, are integers that represent the number of tooths for four gears on a gearwheel.

The computations for all algorithms have been done for 30 independent trials and 1000 function evaluations. The results obtained with LX-BBSCA, LX-BBO, BBO, and SCA, as well as other metaheuristic algorithms, are reported in Table 5. The results show the superiority of the LX-BBSCA algorithm in solving the gear train problem by obtaining a better fitness value. Figure 1 shows the convergence graph for the problem.

**Table 5.** Comparison of different algorithms for the gear train problem.

| Technique | $x_1$ | $x_2$ | $x_3$ | $x_4$ | Objective function value | The largest value of the objective function |
|---|---|---|---|---|---|---|
| LX-BBSCA | 12 | 60 | 17.611 | 24.41 | $1.507 * 10^{-11}$ | 1000 |
| LX-BBO | 23 | 39 | 12 | 50.08 | $2.58 * 10^{-11}$ | 1000 |
| BBO | 14.66 | 60 | 16.67 | 28.25 | $2.68 * 10^{-11}$ | 1000 |
| SCA | 18 | 26 | 12 | 60 | $7.5*10^{-10}$ | 300 |
| GA | 17 | 33 | 14 | 50 | $1.3*10^{-09}$ | NA |
| CS | 16 | 43 | 19 | 49 | $2.7*10^{-12}$ | 5,000 |
| ABC | 16 | 44 | 19 | 49 | $1*10^{-05}$ | 40,000 |
| GWO | 17 | 54 | 22 | 48 | $1.1*10^{-10}$ | 700 |

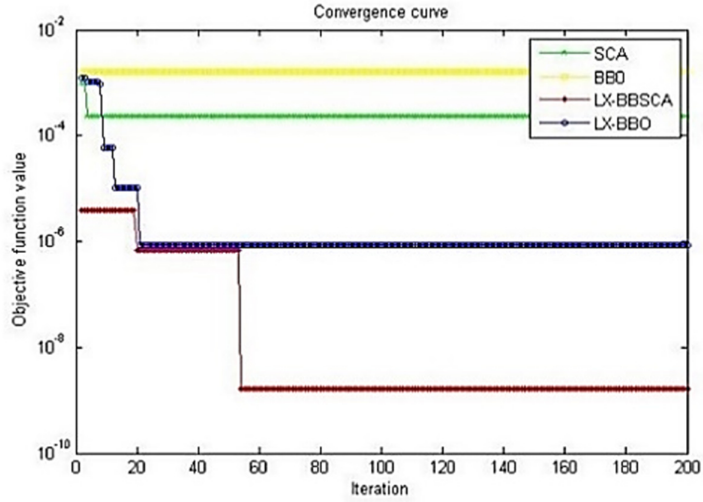

**Figure 1.** Convergence graph of different algorithms for the gear train problem.

*5.2. Gas production problem*

The level of gas production capacity that maximizes the effectiveness and profitability of a gas production facility or operation is referred to as the optimal capacity of gas production. Finding the right balance between production costs, market demand, infrastructural capabilities, and other elements that have an impact on the efficiency of the gas production process is necessary to determine the appropriate capacity. This problem consists in choosing the best production capacities for a system that generates and contains gas. The problem can be expressed as follows:

$$\text{minimize } f_2(x) = 61.8 + 5.72x_1 + 0.2623[(40 - x_1)ln\left(\frac{x_2}{200}\right)]^{-0.85}$$

$$+0.087(40 - x_1)ln\left(\frac{x_2}{200}\right) + 700.23x_2^{-0.75} \quad (9)$$

$$\text{s. t.} \qquad 17.5 \leq x_1 \leq 40$$

$$300 \leq x_2 \leq 600 \qquad (10)$$

where $x_1$ and $x_2$, are the capacity of production facilities.

The results obtained with different algorithms are reported in Table 6. The obtained results demonstrate that the LX-BBSCA outperforms other techniques by obtaining a higher objective function value. Figure 2 shows the convergence graph for the problem.

**Table 6.** Comparison of different algorithms for the gas production problem.

| Technique | $x_1$ | $x_2$ | Objective function value | The largest value of the objective function |
| --- | --- | --- | --- | --- |
| LX-BBSCA | 26 | 600 | 170.1 | 1000 |
| LX-BBO | 17.5 | 400 | 166.3 | 1000 |
| BBO | 40 | 529 | 168.8 | 1000 |
| SCA | 17.5 | 600 | 169.843 | 300 |
| PSO | 17.5 | 600 | 169.843 | 342 |
| ABC | 17.5 | 600 | 169.843 | 319 |
| DE | 17.5 | 593 | 169.996 | 324 |
| GWO | 17.51 | 600 | 169.818 | 342 |

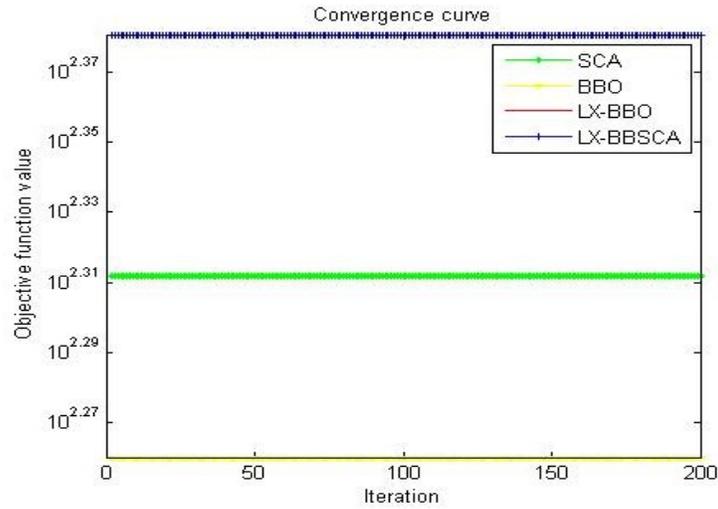

**Figure 2.** Convergence graph of different algorithms for the gas production problem.

*5.3. Beam design problem*

When designing a beam, the right dimensions, materials, and configurations must be chosen to guarantee that the beam can safely carry the applied loads and satisfy certain design requirements. Beams are structural components used to distribute and carry loads to supports. The problem of beam design is multi-objective and can be stated as follows:

1. The area of the beam's cross-section that for a given length reflects its volume.
2. The beam's static deflection for the displacement under force P.

Both criteria should be kept to a minimum. These criteria are obviously incompatible because the smallest cross-section area yields the greatest deflection and vice-versa. The mathematical model of this problem is as follows:

$$\text{minimize } f_3(x) = \frac{5000}{\frac{x_2(x_2-2x_4)}{12} + \frac{x_1 x_4^3}{6} + 2x_1 x_1 \left(\frac{x_2-x_4}{2}\right)^2}, \quad x = (x_1, x_2, x_3, x_4) \tag{11}$$

$$\text{s.t.} \quad g_1(x) = 2x_1 x_3 + x_3(x_2 - 2x_4) \leq 300$$

$$g_2(x) = \frac{18 x_2 \times 10^4}{x_3(x_2 - 2x_4)^3 + 2x_1 x_3(4x_4^2 + 3x_2(x_2 - 2x_4))}$$

$$+ \frac{15 x_2 10^3}{(x_2 - 2x_4)x_3^3 + 2x_1 x_1^3} \leq 6$$

$$10 \leq x_1 \leq 80$$

$$10 \leq x_2 \leq 50 \tag{12}$$

$$0.9 \leq x_3, x_4 \leq 5$$

where $x_1$, $x_2$, $x_3$, and $x_4$, are the height or width of different beams.

The LX-BBSCA, LX-BBO, SCA, and BBO, as well as other metaheuristics, are used to solve the problem with 1000 evaluation functions. The results for this problem are given in Table 7. Figure 3 shows the convergence graph for the problem.

**Table 7.** Comparison of different algorithms for the beam design problem.

| Technique | $x_1$ | $x_2$ | $x_3$ | $x_4$ | Objective function value |
|---|---|---|---|---|---|
| LX-BBSCA | 80 | 50 | 0.9 | 0.5 | 0.0123 |
| LX-BBO | 80 | 50 | 0.9 | 0.5 | 0.0123 |
| BBO | 80 | 50 | 0.9 | 0.5 | 0.0123 |
| SCA | 80 | 50 | 0.9 | 0.5 | 0.0123 |
| CS | 80 | 50 | 0.9 | 2.32 | 0.013 |
| GWO | 80 | 50 | 0.9 | 2.32 | 0.013 |

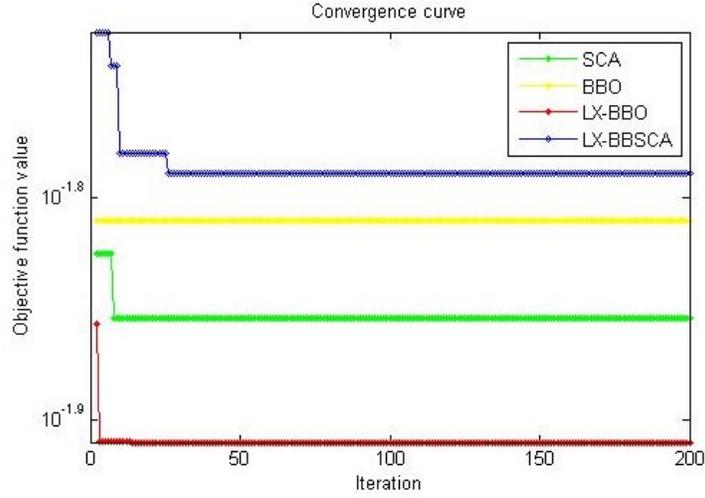

**Figure 3.** Convergence graph of different algorithms for the beam design problem.

*5.4. Space capsule problem*

The enormous engineering task of designing a life support system for a space capsule entails assuring the survival and well-being of astronauts throughout their mission in the hostile environment of space. Everything an astronaut needs to breathe, eat, drink, maintain body temperature, eliminate waste, and stay safe must be provided by the life support system of the space capsule. Due to the system's redundant component's entirely recidivist structure, this dependability issue is challenging to eradicate. This system is made up of four modules, each with a reliability of $r_j$, $j = 1,2,3,4$. It necessitates the success of a single path, and as such, it is split into two subsystems. Each subsystem comprises components arranged in series with component two, creating two identical paths. In this configuration, component three acts as a backup for the third path within this pair, while component one is complemented by parallel component four. The phase formed by elements one and four is in series with element two. Finally, these two identical paths run in parallel, ensuring that one of them guarantees the output. As a final expression, the reliability of the overall system can be defined as:

$$R_5 = 1 - r_3 \times [(1 - r_1)(1 - r_4)]^2 - (1 - r_3)[1 - r_2\{1 - (1 - r_1)(1 - r_4)\}]^2, \quad (13)$$

while the system cost can be calculated as:

$$C_5 = 2k_1 r_1^{p1} + 2k_2 r_2^{p2} + 2k_3 r_3^{p3} + 2k_4 r_4^{p4}, \quad (14)$$

where $k_1 = k_2 = 100$, $k_3 = 200$, $k_4 = 150$, and $p_1 = p_2 = p_3 = p_4 = 0.6$.

The problem can be mathematically expressed as follows:

$$\text{minimize } f_4 = C_5 = 2k_1 r_1^{p1} + 2k_2 r_2^{p2} + 2k_3 r_3^{p3} + 2k_4 r_4^{p4} \quad (15)$$

$$\text{s.t.} \quad 0.5 \leq r_j \leq 1 (\forall j = 1,2,3,4)$$

$$0.9 \leq R_5 \leq 1 \qquad (16)$$

The results of the different algorithms are given in Table 8. Also, Figure 4 shows the convergence graph for the problem.

**Table 8.** Comparison of different algorithms for the space capsule problem.

| Technique | $r_1$ | $r_2$ | $r_3$ | $r_4$ | Cost | Reliability |
|---|---|---|---|---|---|---|
| LX-BBSCA | 0.5 | 0.5 | 0.5 | 0.5 | 725.09 | 0.9991 |
| LX-BBO | 0.5 | 0.5 | 0.5 | 0.5 | 725.08 | 0.9871 |
| BBO | 0.5 | 0.5 | 0.5 | 0.5 | 725.08 | 0.9645 |
| SCA | 0.5 | 0.5 | 0.5 | 0.5 | 725.08 | 0.9528 |
| SA | 0.5009 | 0.8377 | 0.5005 | 0.5001 | 641.903 | 0.9001 |
| GWO | 0.5001 | 0.8389 | 0.5001 | 0.5001 | 641.831 | 0.9 |

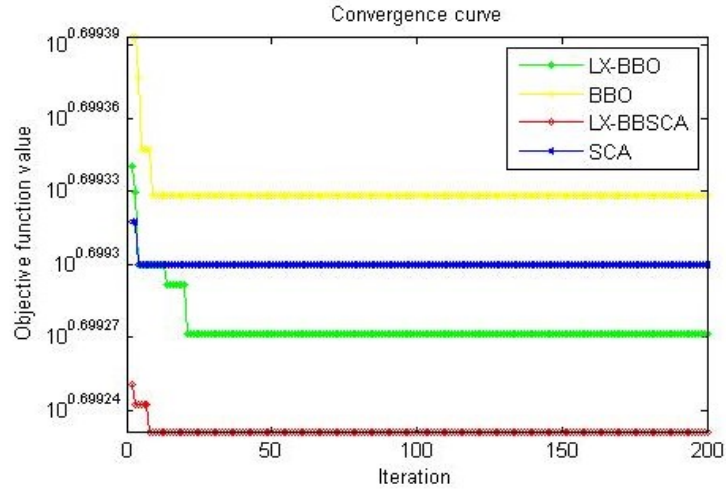

**Figure 4.** Convergence graph of different algorithms for the space capsule problem.

*5.5. Complex bridge network problem*

An interconnected network of bridges that can efficiently handle traffic flow, address geographic constraints, and satisfy numerous engineering requirements is required when designing a complex bridge network. The current issue is a nonlinear multi-objective optimization problem, and the ideal solution is one that maximizes reliability while minimizing expenses. The system consists of five dimensions, each with a reliability of $r_j$, $j = 1,2,\ldots,5$.

$$\text{minimize } f_5 = C_5 = \sum_{j=1}^{5} a_j \, exp\left(\frac{\beta_j}{1-r_j}\right) \tag{17}$$

$$s.t. \quad 0.5 \leq r_j \leq 1 \quad \forall j = 1,2,3,4$$

$$0.9 \leq r_5 \leq 1$$

$$a_j = 1 \text{ and } \beta_j = 0.0003 \quad \forall j = 1,2,3,4,5$$

$$\begin{aligned} R_S(r_1,r_2,r_3,r_4,r_5) &= r_1r_4 + r_2r_5 + r_1r_3r_5 + r_2r_3r_4 \\ &\quad -r_2r_3r_4r_5 + 2r_1r_2r_3r_4r_5 - r_1r_3r_4r_5 \\ &\quad -r_1r_2r_3r_5 - r_1r_2r_4r_5 - r_1r_2r_3r_4 \end{aligned} \tag{18}$$

Table 9 presents the results for the algorithms LX-BBSCA, LX-BBO, BBO, and SCA. Figure 5 depicts the convergence graph for the problem.

**Table 9.** Comparison of different algorithms for the complex bridge network problem.

| Technique | $r_1$ | $r_2$ | $r_3$ | $r_4$ | $r_5$ | Cost | Reliability |
|---|---|---|---|---|---|---|---|
| LX-BBSCA | 0.5 | 0.5 | 0.5 | 0.5 | 0.50 | 5.003 | 0.9901 |
| LX-BBO | 0.5 | 0.5 | 0.5 | 0.5 | 0.50 | 5.003 | 0.9540 |
| BBO | 0.5 | 0.5 | 0.5 | 0.5 | 0.50 | 5.003 | 0.9566 |
| SCA | 0.5 | 0.5 | 0.5 | 0.5 | 0.50 | 5.003 | 0.9508 |
| RST | 0.9392 | 0.9345 | 0.7715 | 0.939 | 0.92 | 5.020 | 0.9900 |
| CSA | 0.9349 | 0.9355 | 0.7853 | 0.941 | 0.92 | 5.019 | 0.9000 |
| GWO | 0.9595 | 0.8798 | 0.7908 | 0.973 | 0.89 | 5.025 | 0.9911 |

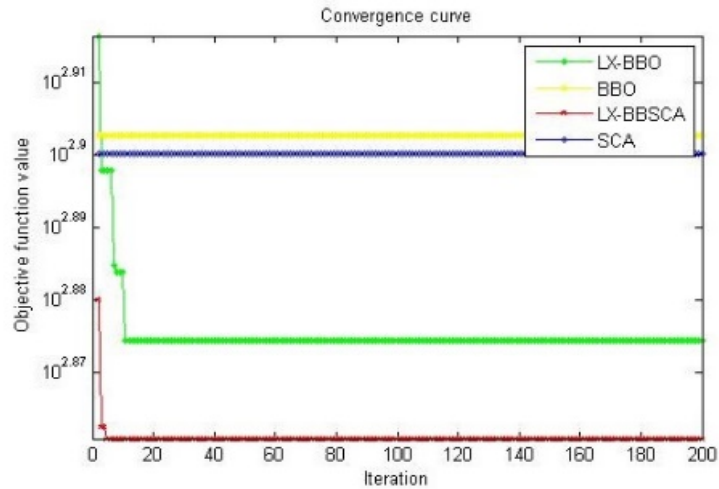

**Figure 5.** Convergence graph of different algorithms for the complex bridge problem.

## 6. Conclusions

In this study, a hybrid metaheuristic-driven optimization algorithm LX-BBSCA has been created by combining SCA and LX-BBO, which coordinates both local and global search and incorporates the advantages of both algorithms. As a result, the novel optimization algorithm has an effective algorithm for searching the solution space. We can draw the conclusion that the new hybrid technique is capable of maintaining a delicate balance between oppression and research based on the convergence performance of the suggested algorithm. Five engineering design issues were resolved using the suggested algorithm. The new hybrid approach surpasses the existing techniques in terms of convergence speed and solution quality. According to the quantitative statistical data, it can be demonstrated that the proposed algorithm may be employed as a successful and proficient ensemble metaheuristic algorithm for real-world tasks involving a complicated search area.